\documentstyle[12pt]{article}
\newtheorem{theorem}{Theorem}
\newtheorem{prop}[theorem]{Proposition}
\newtheorem{lem}[theorem]{Lemma}

\def\grad{\mbox{\rm grad}}
\def\diag{\mbox{\rm diag}}
\def\mod{\mbox{\rm\ mod }}
\def\Blat{\mbox{\it \raise2pt\hbox{"}\kern-2pt H}}
\def\lvup{\rlap{\ ${}^{\ell\atop{\hbox{${}^{\vee}$}}}$}\cdots}
\begin{document}
\center{\Large{\bf Hyperbolic Cauchy problem and Leray's residue
formula.} }
\vspace{1pc}  \center{\large{ Susumu Tanab\'e}}   %
\vspace{2pc}
\center{\it Dedicated to the memory of Bogdan Ziemian}

\abstract {
We give an algebraic  descriptions of (wave) fronts  that  appear  in
strictly hyperbolic Cauchy problem.  Concrete
form of defining function of  wave front issued
from initial algebraic variety is obtained  by the aid of Gauss-Manin
systems satisfied by Leray's residues.
}
\section* {0. Introduction.}
In one of his last works \cite{Ziem}, Prof. Bogdan Ziemian pursued  a
possibility to express fundamental solutions to PDE  by  the  aid  of
Leray residues.  He used this technique  to  write  down  the  Mellin
transform of fundamental solutions to Fuchsian type  PDE  and  proved
that these solutions belong to  so  called  a  class  of  generalized
analytic functions(GAF).
\par In
this note, we show that the advantage to make use of Leray's  residue
formula in hyperbolic Cauchy problem is the fact that it  facilitates
the
calculus of the representative integrals (i.e. a basis of certain
cohomology group) whose summation gives fundamental solutions
to the Cauchy problem.
\par   Let us explain in short how the Leray's  residue  formula  can
be applied to the constrution of fundamental solutions. Let us denote
$ V_x =\{\xi \in {\bf C}^n; F(\xi,x')+x_0=0 \}$ a complex variety of
dimension $n-1$  depending on $x =(x',x_0)$ $\in {\bf R}^m$   defined
by  a  polynomial  $F(\xi,x')+x_0$.  Choose  a continuous family   of    
cycles  $\gamma_x$$\in
H_{n-1}(V_x).$  For $x$ "in generic position" the variety $V_x$ is
smooth and $\gamma_x$ is a family of cycles depending continuously on
$x.$
Suppose that $a(\xi)$ is a
smooth function defined on ${\bf C}^n.$ In this situation the following
equality is called  the Leray's residue formula:  $$   I_a(x)   =
\int_{\gamma_x} a(\xi) \frac{d\xi}{dF} = \frac 1{2\pi i}\int_{\partial
\gamma_x} \frac{a(\xi)}{F(\xi,x')+x_0} d\xi, \leqno(0.1)$$ where
$\partial \gamma_x$ $\in H_n({\bf C}^n \setminus V_x)$ so called Leray
coboundary of  the  cycle  $\gamma_x$ which  is  homologically
equivalent to  a $S^1-$   bundle  over $\gamma_x.$  See  \cite{Ph},
\cite{Vas}.
If $a(\xi)d\xi = d\psi \wedge  dF$  for  some  $\psi  \in
\Omega^{n-2}_{\bf C^n},$ then the  integral  $I_a(x)$  defined  above
must be constantly zero:
$$ I_a(x)
=   \int_{\gamma_x}   \frac{d\psi \wedge dF}{dF}    =
\int_{ \gamma_x} d\psi=0,
\leqno(0.2)$$
that one can see by Stoke's theorem. 

\vspace{1pc}

\noindent\hskip2in\hrule\vskip0.1in
\footnotesize{AMS Subject Classification:  Primary 35L25, 58G17, 33C75;
Secondary 32S40, 78A05, 33C20.  Partially supported by Max Planck Institut
f\"ur Mathematik,Bonn.  }
\normalsize

\newpage
\noindent
Evidently
$$  \int_{\gamma_x } F(\xi,x') a(\xi)    \frac{d\xi}{dF}    =
-x_0  \int_{\gamma_x} a(\xi) \frac{d\xi}{dF}.
\leqno(0.3)$$
After  (0.2)  and  (0.3),  we  conclude  that  the  important   forms
$a(\xi)d\xi$ that will give non zero $I_a(x)$ not expressed by other
Leray's residue must be of the space
$$ \Blat = \frac{\Omega^{n}_{{\bf C}^n}}{dF\wedge
d\Omega^{n-2}_{{\bf C}^n} + F\Omega^{n}_{{\bf C}^n}}.\leqno(0.4)$$

Furthermore (0.1) yields the relation
$$  \int_{\gamma_x}   \frac{d\omega}{dF}    =
\frac{d}{dx_0}\int_{ \gamma_x} \omega,
\leqno(0.5)$$ for $\omega \in \Omega^{n-1}_{{\bf  C}^n}.$
That  is  to
say a differential equation satisfied by $I_a(x)$ does not depend  on
the choice of cycle along which one defines the integral.

In  his  famous  work  \cite{Br},  E.Brieskorn  has  shown  that  for
$F(\xi,x')$ whose singular fibre $V_0$ defines  a hypersurface isolated
singularity, the space \Blat\  is a vector space of finite
dimension $\mu$ that coincides  with  the
Milnor number of the singularity $ V_0 =\{\xi \in {\bf C}^n;
F(\xi,0)=0 \}.$
Thus if  a fundamental solution $I_a(x)$ is expressed by a sum of
integrals of
forms  $b_1  d\xi,  \cdots,b_\mu  d\xi$ $\in \Blat$,  for  certain
$F(\xi,x')$ we can expect that analytic properties of the  fundamental
solution  can  be  deduced  from   those  of  $I_{b_1}(x),    \cdots,
I_{b_\mu}(x)$ which in their turn can be described by informations of
the singularity $V_0.$

To pursue further
this study,  we propose to make use of the  Gauss-Manin
system associated with fibre bundle structure that naturally arises in
integration. Our main tool is concrete expressions of the
overdetermined differential systems obtained from non trivial relations
between base elements of \Blat\ for complete intersection
singularities  (Proposition~\ref{prop21}, Theorem~\ref{theorem21}).
In our situation, Leray's residue formula can be written down like
(2.10) below. The main Theorem ~\ref{theorem22} directly follows
the Theorem~\ref{theorem21} on reducing the situation to
a specific mapping (2.15).
\par   In our former work \cite {T4}
we illustrated in concrete examples the possibility to interprete
the fundamental solution to the Cauchy problem associated with
the wave operator $$ P(D_t,D_x)= (\frac{\partial}{\partial t})^2 -
\sum_{j=1}^n( \frac {\partial}{ \partial x_j})^2,$$
as generalized hypergeometric
functions. In this note we consider  the  Cauchy  problem  associated
with strictly hyperbolic operators with constant coefficients
in general. This
procedure has been supported by the general theory of the
Gauss-Manin systems for isolated complete intersection  singularities
\cite {Gr}. More systematic explanation of this situation from singularity
theoretical  point  of  view  is  given  in  \cite{AT2},   \cite{T2},
\cite{T3} and \cite{T4}.

\section* {1. Preliminaries on Cauchy problem.}  In this section we
prepare fundamental notations and lemmata to
 develop our studies in further sections.
Let $P(D_t,D_x)$ be a strictly  hyperbolic  operator  with
constant coefficients of degree $m$ i.e. its total symbol
$$ P(\tau,\xi) = \tau^{m}  +  \sum_{i=1}^{m}  P_{m-i}(\xi)\tau^{m-i},
$$
$$ P_{m-i}(\xi)= \sum_{|\alpha|= i} P_{m-i,\alpha}\xi^\alpha \in {\bf
R}[\xi],$$
satisfies the following decomposition
$$  P(\tau,\xi) = \prod_{j=1}^{m}(\tau - \lambda_j(\xi)), $$
such that $\lambda_j(\xi) \in {\bf R}$ if $\xi \in {\bf R}^n$ and
$\lambda_i(\xi) \not =$ $\lambda_j(\xi)$ for $i \not= j,$ $\xi
\in{\bf R}^n \setminus \{0\}.$

Without loss of generality, we suppose that $P(\tau,\xi)$ is an irreducible
polynomial of  ${\bf R}[\tau,\xi].$ If $P(\tau,\xi) = P_1(\tau,\xi)
\cdot P_2(\tau,\xi)$
fundamental solutions of
$P(D_t,D_x)$  is   a   sum   of   those   of   $P_1(D_t, D_x)$   and
$P_2(D_t, D_x)$ provided that their characteristic roots are mutually
distinct out of the  origin.  Let us consider the following Cauchy
problem (C.P.):  \begin{flushleft} $ {\rm (C.P.)}  \left \{
 \begin{array}{ccc} P(D_t,D_x)u(t,x) &= &0  \\ D_t^{m-1} u(0,x) & =& v
(x)  \\ D_t^{m-j}u(0,x)&=&  0 , 2 \leq j \leq m  \end{array} \right. $
\end{flushleft} where $t\in {\bf R},  x \in {\bf
R}^n,  D_t = \frac {\partial}{i \partial t}, D_x=
( \frac{\partial}{i \partial x_1}, \cdots,
\frac{\partial}{i \partial x_n}), i= \surd -1.$
We will study the Cauchy problem
(C.P.) under the following conditions  (C.1), (C.2), (C.3) imposed on the
initial data.  In order to  describe these conditions, we use
notation $\chi^{\epsilon}_q(z)\;\;(\epsilon = \pm 1)$ which
stands for the following distributions defined as boundary values of
an  analytic function on ${\bf C}_z^1 \setminus \{0\}$ (cf.
\cite{G}) :$$\chi^{\epsilon}_q(z)= \chi_q(z + i 0) + \epsilon \chi_q(z
- i 0),$$ where $$  \chi_q(z  )  =\Gamma(-q) z^q \;\; q \notin  {\bf
Z}\;\;{\rm or}\;q\;\mbox{negative integer}$$ $$ \chi_q(z ) = \frac
{z^q}{q!}(-\log z +  C_q) \;\;\;\;q{\;\mbox{positive integer}}$$
where $ C_0 = 0 , \; C_q = C_{q-1} + 1/q.$ Remark that $$ \frac{d}{dz}
\chi_q(z) = \chi_{q-1}(z).  $$ (C.1) {\it The initial data are given by
a distribution of finite order with singular support (see Definition
2.2.3 of \cite{H}) located on cotangent bundle of a smooth algebraic
surface $S := \{ x \in {\bf R}^n : F(x) -s =0 \}$ defined by a real
polynomial $F(x),$ $$ v(x) = g(x)\chi^{\epsilon}_q(F(x) -s)  $$ \it
with a smooth function $g(x)$.} \par Further we shall denote the
singular support of a distribution $v(x)$ by $S.S. v(x).$ \par We
impose several technical conditions also. These conditions will be used
so that the reasoning on the isolated complete  intersection
singularities can be applied to our (C.P.) \par
(C.2)(Quasihomogeneity){\it  There exists a set of positive integers
$(w_1, \cdots, w_n),$ that satisfies 1) $w_i \not = w_j$ for certain
pair $1 \leq i \not = j  \leq n, $  2) $G.C.D.(w_1, \cdots, w_n)=1$ and
3) for a positive integer  $w(F),$ $$ (\sum_{1\leq j \leq n} w_j
 x_j\frac {\partial }{\partial x_j}) F(x) = w(F)F(x),$$ holds.} \par
The above contdition 1) plays an essential r\^ole in establishing the
Lemma 8 below.

\par (C.3){\it The following is a vector space of finite
dimension:} $$ {\bf R}[x]/{\cal I}$$ {\it where ${\cal I}= \langle
F(x), \frac {\partial F}{\partial x_1}(x), \ldots, \frac {\partial
F}{\partial x_n}(x) \rangle$ ( ideal generated by  the entries).}  Let
us introduce the following notations.  \par $a)$ The phase function
$\psi (x,t,z)$ is defined as follows:  $$\psi (x,t,z) =  P( \langle
x-z, \grad_z F(z) \rangle, t \; \grad_z F(z) ).  $$ \par $b)$ The paired
oscillatory integrals studied in \cite{G} defined for the phase
function $\psi (x,t,z)$ introduced in  a): $$ I^{\epsilon}_{p}(x,t,s) =
\int_{\{ F(z) = s \}} H_{ p}(z)  \chi_p^{\epsilon}(\psi(x,t,z))
\frac{dz}{dF}, \;\;$$ with regular amplitude functions of
pseudo-differential operator $$ H_{p}(z) \sim \sum_{r=p-m}^{-\infty}
 h_{ p,r}(z) \in S^{p-m}({\bf R}^n),$$ in which  $h_{p,r}(z)$ is
homogeneous of order $r$ for large values of $z.$ One shall understand
$I^{\epsilon}_{p}(x,t,s)$ as the Gel'fand-Leray integral  ( see 1.5
\cite{Br}) defined on the real algebraic set  $S=\{  z\in {\bf R}^n ;
F(z) =s\}. $ \par $c)$ The function $\phi(x,t,s) $   denotes a defining
function of the (wave) front $\Sigma$ issued from $S$ determined by
$(C.P.).$

\begin{prop}      The following assertions hold for  solutions to
the Cauchy problem  (C.P.) with the notations introduced in a),b) and c) as
above. \par d) The solution $u(x,t)$ to the Cauchy problem (C.P.)
admits  an
asymptotic expansion  $$  u(x,t) \sim \sum_{j=0}^{\infty}
I^{\epsilon}_{-\frac{n}{2} + q+ j}(x,t,s). \leqno(1.1)$$ That is to say, for every
$N>>0$ there exists $C_N >0 $ such that $$ \mid u(x,t) - \sum_{j=0}^{N}
 I^{\epsilon}_{-\frac {n} {2} + q+ j} (x,t,s) \mid \leq C_N \mid \phi(x,t,s)
 \mid^{q+N+1}, \leqno(1.2) $$ in the neighbourhood of
$S.S.u(x,t).$  \label{prop11} \end{prop}  {\bf Proof} We give only a sketch of
proof while a detailed one shall appear in \cite{T3}, \cite{T4}. First of
all we show
that the phase function of the integrals  $I^{\epsilon}_{p}(x,t,s)$
coming into
$a)$ is given by $c).$  In solving the Hamilton-Jacobi equation associated with
the Hamiltonian  $\tau  -  \lambda_\kappa  (\xi),$  $1  \leq  \kappa
\leq m$
(in a symplectic coordinate with canonical
symplectic form $dt\wedge d\tau + \sum_{j=1}^n dx_j \wedge d\xi_j$. For the
symplectic  geometry see Chapter XXI \cite{H}),  \begin{center} $\left
\{ \begin{array}{ccc} \dot t& = & 1 \\ \dot \tau& = & 0 \\ \dot x_i& = &
\frac { \partial \lambda_\kappa (\xi) }{\partial \xi _i} \\ \dot \xi_i&
= &  0  \\
x_i (0) & = & z_i,\;\;\{z\in{\bf R}^n; F(z) =s\}, 1 \leq i \leq n,
\\ \end{array} \right. $ \end{center}   we
get
\begin{flushleft}   $(1.3)\;\;\;\;   x_i   =    t   \frac   {\partial
\lambda_\kappa(\xi)}{\partial \xi _i} + z_i
\;\;{\rm with }  \;\;z \in S=\{  z\in
{\bf R}^n ; F(z) =s\}. $ \end{flushleft} This means that the singularities of
the solutions to (C.P.) lies on the rays (1.3). These lines are interpreted as
rays issued from the initial front $S = \{ z\in {\bf R}^n: F(z)= s  \}$
in
directions determined   by  the  Hamiltonian
$$ P(\tau,\xi) = \tau^{m}  +  \sum_{i=1}^{m}  P_{m-i}(\xi)\tau^{m-i}
= \prod_{j=1}^{m}(\tau - \lambda_j(\xi))$$  in question.

 Consequently
they are expressed by integrals with phase
$$\psi (x,t,z) =
 \prod_{j=1}^{m}\left(\langle x-z, \grad_z F(z) \rangle  -
t\; \lambda_j(\grad_z  F(z))\right)\leqno(1.4)$$  $$ =  P(\langle  x-z,
\grad_z F(z) \rangle,  t\;  \grad_z  F(z)). $$  This  is  the
``minimal'' algebraic  equation  describing  the  wave  front  in  view
of   the irreducibility of the polynomial $P(\tau,\xi).$ Here we remark
that for every $p \in \bf Q$ and $H(z) \in  {\cal D}'({\bf R}^n_z),$ we
have $$ P(D_t, D_x)\int_{S} H(z)(\psi(x,t,z))^p \frac{dz}{dF} = 0.$$
One can prove this equality with the aid of Gauss-Stokes' theorem. Thus
the question is how to find a series of integrals  $$
I^{\epsilon}_{p}(x,t,s) = \int_{S} H_{p}(z) \chi_p^{\epsilon}(\psi
 (x,t,z)) \frac{dz}{dF}, \;\;\; p \in {\bf Q}  $$ whose suitably
converging sum produces a distribution $u(x,t)$ satisfying (C.1).  The
 possibility of an asymptotic expansion (1.2) consisting of terms like
$b)$ can be proven by well known estimates on the stationary phase
 (\cite{H}, Theorem 7.7.12). More precisely, let us remind the  following
lemma.   \begin {lem} Let $(G)^0(x,t,s)$ be  a residue of a smooth function
$G(x,t,z)$ after division by an Jacobi ideal  generated by  $\frac{\partial
\psi(x,t,z)}{\partial z_j},$ $\frac{\partial F(z)}{\partial z_j}, \; 1\leq j
\leq n,$ and $F(z)-s,$ i.e.  $$\leqno(1.5) G(x,t,z)= G^0(x,t,s) + \sum_{j=1}^n
f_j(x,t,s,z) \frac{\partial \psi(x,t,z)}{\partial z_j} +$$ $$ + \sum_{j=1}^n
g_j(x,t,s,z) \frac{\partial F(z)}{\partial z_j}  + h(x,t,s,z)(F(z)
-s),    $$    with    some    smooth     functions      $h(x,t,s,z),$
$f_j(x,t,s,z),$ $g_j (x, t, s, z ),$ $1\leq j \leq n.$
\par Then
for every smooth function $a(z)$  the following
asymptotic estimate with some $C_N >0$ holds  in the neighbourhood of the
 wave front  $\Sigma =\{ (x,t) \in {\bf R}^{n+1};$  $(\psi)^0 (x,t,s) = 0
\}$:
\begin{flushleft} $(1.6) \mid \int_S  a(z) \chi^{\epsilon}_q(\psi(x,t,z)) dz -
\sum_{j= 0}^N (L_{\psi,j}a)^0(x,t,s) \chi^{\epsilon}_{q + n/2 + j} ((\psi)^0
(x,t,s))\mid$ \end{flushleft}
$$ < C_N \mid (\psi)^0 (x,t,s) \mid^{n/2 + N +1+q} $$ with differential
operators $ L_{\psi,j}$ of degree $2j.$ Furthermore
 we have: $$(L_{\psi,0}a)^0(x,t,s)  = i^{\frac{n}{2}}(2\pi)^{\frac{n-1}{2}}
(a)^o(x,t,s) \mid \det  (\frac{\psi_{zz}}{2 \pi
i})^0(x,t,s)\mid^{-1/2} .$$ \label{lem11} \end{lem}
 In the literature concerning the singularity theory, one often calls
the correspondence $G(x,t,z)$ $\rightarrow$ $G^0(x,t,s)$
Lyashko-Loojenga mapping.  \par Let us briefly  sketch proof of the
lemma. Malgrange's division theorem yields the decomposition (1.6) in
connexion with the fact that the follwing ${\cal O}_{{\bf C}^n}-$
module is a finite dimensional vector space under the assumption (C.3):
$$ \frac{\Omega^n_{{\bf C}^n}}{dF(z) \wedge\Omega^{n-1}_{{\bf C}^n}+
d\psi(0,0,z) \wedge\Omega^{n-1}_{{\bf C}^n}+ d_z\Omega^{n-1}_{{\bf
C}^n} + F(z) \wedge\Omega^{n}_{{\bf C}^n} }.$$ Further it suffices to
apply above mentioned  stationary phase method. \par  After
Lemma~\ref{lem11} the function $ (\psi)^0(x,t,s)$ can be given by
(1.5) for $G(x,t,z)$ $ = $ $\psi(x,t,z)$ in (1.4). As a definig function
$\phi(x,t,s)$ of the wave  front issued from
$S,$   one  shall  take  a  polynomial  such  that  $\{   (x,t)   \in
{\bf R}^{n+1}; (\psi)^0(x,t,s) =0 \}$  $\subset
\{   (x,t)   \in {\bf R}^{n+1}; \phi(x,t,s) =0 \}$ and $S=
\{x:\phi(x,0,s)=0 \}.$ We adopt such $\phi (x,t,s)  $ as needed one in
$c).$ \par It remains to justify asymptotic estimates in $b)$ and $d).$
This can be achieved in view of (1.6) and well known construction of an
elementary solution to  strictly hyperbolic Cauchy problem (see for
example \cite{Ham}, \cite{HLW}). Especially due to the choice of
$\phi(x,t,s),$ the inequality (1.2) satisfied with $\psi^0(x,t,s)$
in a local context holds with $\phi(x,t,s).$ Hence the
assertion follows.  {\bf Q.E.D.} \par We formulate a simple lemma
before introducing necessary notations.  \begin{lem} Under the
assumptions (C.2), (C.3) imposed on $F(x)$ there exists a collection of
polynomials degree $\leq m,$ $W_1(x,t),$ $\cdots,$ $ W_{\mu'}(x,t),$
with  $\mu'$ an integer smaller than $m^n \prod_{i=1}^n$
$(\frac{w(F)}{w_i}),$ satisfying
$$
\psi(x,t,z) = { \langle z, \grad_z
F(z)) \rangle }^m  + \sum_{i=1}^{\mu'} W_i(x,t)z^{\alpha^{(i)}}
\leqno(1.7)
$$
for $\psi(x,t,z)$ of (1.4).  Here $\alpha^{(i)}  =
(\alpha^{(i)}_1, \cdots, \alpha^{(i)}_n  )  \in ({\bf Z}_{\geq 0})^n$
stands for  multi-index under restriction $\sum_{j=1}^n w_j
\alpha^{(i)}_j < m \cdot w(F).$ \label{lem12} \end{lem} The proof
follows direct calculation of (1.4).  Quasihomogeneous type of $F(z)$
yields the estimate on term number $\mu'.$ Let us denote by  $$w
(z^{\alpha^{(i)}})= \sum_{j=1}^n w_j \alpha^{(i)}_j,$$ quasihomogeneous
weight of monomial $z^{\alpha^{(i)}}$ for $\alpha^{(i)} \in {\bf N}^n.$
In terms of the quasihomogeneous weight we distinguish two cases.  \par
{\bf Case 1} If there is a term with $w(z^{\alpha^{(i)}}) =0$ in (1.7), let us
mark it as $\alpha^{(1)}$ and define a polynomial $$ f_1(y(x,t),z) = {
\langle z, \grad_z F(z)) \rangle }^m  + \sum_{i=1}^{\mu'}
y_i(x,t)z^{\alpha^{(i)}}.$$ Here $y_i(x,t)=W_i(x,t), 1 \leq i \leq
\mu',$ for polynomials introduced in Lemma~\ref{lem12}, (1.7).  \par
{\bf Case 2} If all terms of (1.7) has positive weight, we shall define
$$ f_1(y(x,t),z)  =  {  \langle  z,  \grad_z   F(z))   \rangle   }^m
+ \sum_{i=2}^{\mu' +1} y_i(x,t)z^{\alpha^{(i-1)}} +y_1.$$ with
$y_{i+1}(x,t)=W_i(x,t), 1 \leq i \leq \mu'.$ \par For the sake of
simplicity we adopt notation $\mu = \mu'$ for Case 1 and $\mu = \mu'
+1$ for Case 2.  \par Further we define integrals $$ I_{p}(y(x,t),s) =
\int_{S} H_{p}(z) \chi_p(f_{1}(y(x,t),z)) \frac{dz}{dF}.  \leqno(1.8)$$
Hence if one denotes $y' = (  y_2(x,t), \cdots, y_\mu(x,t) ),$ $$
I^{\epsilon}_{p}(x,t,s) = I_{p}(y_1 +i0, y'(x,t),s) + \epsilon
I_{p}(y_1 -i0, y'(x,t),s) $$ on understanding that the boundary value
shall be taken at $y_1 =0$ in the above mentioned Case 2. Thus it is essential
to study  $I_{p}(y(x,t),s)$ of (1.8) to estimate asymptotic behaviour of  $
I^{\epsilon}_{p}(x,t,s). $ From now on we shall regard the integral (1.8) as a
function in variables $y(x,t) = (y_1, y_2(x,t), \cdots, y_\mu(x,t)).$
 Therefore our main concern will be to investigate the differential equations
 that satisfy $I_{p}(y,s)$ corresponding to various amplitudes $H_p(z)$ 
with the aid  of the Gauss-Manin connexions associated  to  
complete intersection
singularities.
\section*{2. Gauss-Manin connexions for
quasihomogeneous  complete intersections.}   Here we propose to study the
integrals $I_{p}(y,s)$ defined in (1.8) by means of the Gauss-Manin system
 associated with complete intersection singularities.  In effect, it is well
known that the Gauss-Manin connexion can be  defined on the relative de Rham
cohomology  groups. Instead of that here we propose to calculate it on spaces
so called Brieskorn lattices (see \cite{Br}, \cite{Gr}).  \par The formulation
of this section is a modification of \cite{T2}, \S 1 adapted to our situation.

 Let us observe a mapping between complex manifolds $ X = ( {\bf C}^{N+K},
0),Y = ( {\bf C}^K, 0),$  $$ f: X \rightarrow Y $$ that defines an isolated
quasihomogeneous complete intersection singularity  at the origin. That is to
say, if we denote  $$X_y :=\{ u \in X ;  f_0(u) = y_0, \ldots, f_{K-1}(u) =
y_{K-1}  \},\leqno(2.1) $$  then $\dim X_y =N \geq 0$ and the critical set of
mapping $ f: X_0 \rightarrow Y $ is isolated in $X_0.$ Further we assume that
polynomials  $f_0(u), \cdots,f_{K-1}(u)$ are quasihomogeneous i.e. there
exists a collection of positive integers $v_1,\cdots, v_{N+K}$ whose
greatest common divisor equals 1 and  $$ (v_1u_1
\frac{\partial}{\partial u_1}+ \cdots  + v_{N+K} u_{N+K}
\frac{\partial}{\partial u_{N+K}}) f_{\ell}(u) = p_{\ell} f_{\ell}(u),
\;\;\; \ell = 0,1,\cdots, K-1,$$ for certain integers $p_0,\cdots,
p_{K-1}.$ We shall call the vector field $$E= \sum_{i=1}^{N+K} v_i u_i
\frac{\partial}{\partial u_i}  , \leqno(2.2)$$ Euler vector field and
$v_1,$ $\cdots,$ $v_{N+K}$  (resp.  $p_0,$ $\cdots,$ $ p_{K-1}$)
 positive weights of variables $u_1,$ $\cdots,$ $u_{N+K}$ (resp.
polynomials $f_0,$ $\cdots,$ $f_{K-1}$) i.e.  $v_1= w(u_1),$ $p_0 =
w(f_0)$ etc. \par In order to calculate the Gauss-Manin connexion for
isolated complete intersection singularity $ X_0,$ we introduce two
 vector spaces $V$ and $F.$ After Greuel-Hamm \cite{GH}, we look at a
space whose dimension as a vector space  over $\bf C$ is known to be
the Minor number $\mu(X_0)$  of singularity $X_0,$ $$ V : =
\frac{\Omega^{N}_X} { df_0\wedge\Omega^{N-1}_X +\cdots +
df_{K-1}\wedge\Omega^{N-1}_X + d\Omega^{N-1}_X  + f_0 \Omega^{N}_X +
\cdots +  f_{K-1} \Omega^{N}_X }.\leqno(2.3) $$   The second one will
later turn out to be isomorphic to $V$ (see  Proposition
~\ref{prop22}), $$F : = \frac{\Omega^{N+1}_X}{df_0\wedge\Omega^{N}_X
 +\cdots + df_{K-1}\wedge\Omega^{N}_X +  i_{E}(\Omega^{N+2}_X )}.
 \leqno(2.4)$$ Here $i_E$ means the inner contraction with Euler field
$E$ defined by (2.2). The third vector space associated with the
singularity $X_0$ is defined as follows $$ \Phi:=
  \frac{\Omega^{N+K}_X}{df_0\wedge \cdots \wedge df_{K-1} \wedge
\Omega^{N}_X +  f_0\Omega^{N+K}_X + \cdots +
f_{K-1}\Omega^{N+K}_X} .\leqno(2.5)$$ Later we define period integrals as
coupling of forms of $V$ or of $\Phi$ with base element of homology groups
$H_N(X_y).$ We remember also the definition of the
Brieskorn lattice \Blat\ from
\cite{Gr},  $$\Blat =\frac{\Omega^{N+K}_X} {df_0 \wedge \cdots \wedge df_{K-1}
\wedge   d\Omega^{N-1}_X },$$ whose rank as ${\cal O}_Y-$ module equals
the Minor number $\mu(X_0)$ of the  singularity $X_0.$  It is easy to
show \begin{lem}  For  $f_0,\cdots, f_{K-1}$ quasihomogeneous polynomials
defining an  isolated complete intersection singularity, $$\Phi \cong
\Blat/(f_0,\cdots,f_{K-1}).$$  Thus $\dim_{\bf C} \Phi = \mu(X_0).$ \label
{lem21} \end{lem}   Let us fix a set of quasihomogeneous $N+1-$ forms
$\tilde{\omega_1},$ $\cdots$ $\tilde{\omega}_{\mu(X_0)}$ whose residue 
class forms a base of $F.$ Similarly we fix a set of $N+K$ 
quasihomogeneous forms $\phi_1(u)du,$ $\cdots$ $\phi_{\mu(X_0)}(u)du$ 
whose residue class gives a base of $\Phi.$

 From definitions (2.4) and  (2.5)
we easily deduce the following.  
\begin{prop}   For each form
$\tilde{\omega_i},$  one has the following decomposition, 
$$ \tilde{\omega_i}\wedge  df_0 \wedge\lvup\wedge df_{K-1}
= \sum_{j=1}^{\mu(X_0)} P^{(\ell)}_{ij} \phi_j(u) du \mod (df_0 \wedge
\cdots \wedge df_{K-1}\wedge d\Omega^{N-1}_X )\leqno(2.6) $$  with
$P^{(\ell)}_{ij}\in {\bf C}[f_0,\cdots, f_{K-1}]$ and $\phi_j(u)du$
as above  for $1 \leq i,j \leq \mu(X_0), 0 \leq \ell \leq K-1$ and   $df_1
\wedge\lvup \wedge df_{K-1} = \bigwedge_{i\,/\kern-5pt=\ell}
df_i.$ \label{prop21} \end{prop}   
As a matter of fact the righthand side of $(2.6)$ can be 
considered as the element of $\Blat.$ 
From \cite{T2} we remember the
 following 

\begin{prop}   Under the situation and definitions as above,
the mapping, $$ i_E : F \rightarrow V $$  induces an isomorphism.
Consequently  $\dim_{\bf C} F = \dim_{\bf C} V = \mu (X_0)$.
\label{prop22}  \end{prop}     In view of  Proposition ~\ref{prop22},
let us choose a base of  $F$ by  $ \tilde{\omega}_i$ (resp. $V$
by $\omega_i$) such
that    $\tilde{\omega_i}$ $= \frac{1}{\ell_i} d \omega_i$
where $\ell_j$ denotes weight of the form ${\omega}_j.$ Remark that
$i_E \tilde {\omega_i}$ $\equiv $ $\frac{1}{\ell_i}{(d i_E + i_E d )
({\omega}_i)}$ $\equiv$ $ \omega_i$
$, 1 \leq i \leq \mu(X_0)$ in $F.$   To make a transition from $(N+K)-$
forms to period integrals, we introduce meromorphic $N-$forms $\psi_i$
satisfying $$ df_0\wedge \cdots \wedge df_{K-1} \wedge \psi_i =
\phi_i(u) du, \;\;\; 1 \leq i \leq \mu(X_0).$$Then we derive the
following relation from Proposition ~\ref{prop21}, \begin{flushleft}
$(2.7)\;\;\; d{\omega}_j
 = \ell_j {\tilde{\omega}_j}$  $\equiv \ell_j
(\sum_{q=1}^{\mu(X_0)} P_{jq}^{(0)} df_0 \wedge \psi_q  + \cdots $ $+
(-1)^{K-1}\sum_{q=1}^{\mu(X_0)} P_{jq}^{(K-1)} df_{K-1} \wedge \psi_q
)$ mod~$((df_0,\cdots, df_{K-1}) d\Omega_X^{N-1}).$\end{flushleft}
See (2.12) below to see that this relation calculates the ``partial
derivative'' $\frac{d\omega}{df_i}.$  Hence, $$ \omega_j \equiv
i_E(\tilde{\omega}_j) \equiv \leqno(2.8)$$

\begin{flushleft} $\sum_{i=0}^{K-1}
(-1)^{i}[ \sum_{q=1}^{\mu(X_0)} P_{jq}^{(i)}p_i f_i \psi_q  -\sum_{
q=1}^{\mu(X_0)} P_{jq}^{(i)} df_i \wedge i_E(\psi_q )]$  mod~$((df_0,\cdots,
df_{K-1}) i_Ed\Omega_X^{N-1}, (f_0,\cdots,
f_{K-1})d\Omega_X^{N-1}).$ \end{flushleft}

As a consequence
 \begin{flushleft} $(2.9)\;\;\;\; d\omega_j \equiv
\sum_{q=1}^{\mu(X_0)} [\sum_{i=0}^{K-1} (-1)^{i}( d (p_i P_{jq}^{(i)}  f_i) -
w(\psi_q) P_{jq}^{(i)} df_i)]\wedge\psi_q  +
\sum_{q=1}^{\mu(X_0)}[\sum_{i=0}^{K-1} (-1)^{i}  p_i P_{jq}^{(i)} f_i]\wedge
 d\psi_q, \;$ mod~$((df_0,\cdots, df_{K-1}) d\Omega_X^{N-1}),$
\end{flushleft}
where $w(\psi_q), $ quasihomogeneous weight of  form $\psi_q.$   The
expression
(2.9) can be simplified if one lets them couple with a vanishing $N-$cycle,
say $\gamma(y)$ and attains non trivial relations between integrals
 $\int_{\gamma(y)}\psi_q, $ instead of those between forms. One defines so
called  period integral  $I_{\phi_q, \gamma(y)}(y)$ taken along a vanishing
cycle $\gamma(y)$ whose ambiguity in homology class  $H_N(X_y, {\bf Z})$ we do
not care for the moment,  \begin{flushleft} $(2.10)\;\;\;\;
I_{\phi_q,\gamma(y)}(y)
: = \int_{\gamma(y)}\psi_q  = (\frac{1}{2\pi i})^{K}
\int_{\partial \gamma(y)}\frac  {df_0 \wedge \cdots \wedge df_{K-1}  \wedge
\psi_q}{(f_0 -y_0)  \cdots (f_{K-1} - y_{K-1})}$ \end{flushleft} $$ =
(\frac{1}{2\pi i})^{K}  \int_{\partial \gamma(y)}\frac  {\phi_q(u) du }{(f_0
-y_0)\cdots (f_{K-1} - y_{K-1})}, $$  where ${\partial \gamma(y)}$
$\in H_{N+K}({\bf C}^{N+K} \setminus \cup_{i=0}^{K-1} \{f_i = y_i\},
{\bf Z})$ is a cycle obtained by the aid of  Leray's coboundary
operator $\partial.$ That is to say, although $\psi_q$ is in general
a meromorphic form with poles along the critical set of the mapping $f,$
$I_{\phi_q,\gamma(y)}(y)$ can be calculated as an integral of a
holomorphic form on  $\partial \gamma(y).$

One
may consult a booklet by F.Pham \cite{Ph}, or a book by V.A.Vasiliev
\cite{Vas} on the coboundary operator.  One understands  (2.10)  the
Leray's residue formula in our situation (2.1).  \par From (2.8) we can
deduce $$\int_{\gamma(y)}{\omega}_j = \sum_{q=1}^{\mu(X_0)}
[\sum_{i=0}^{K-1} (-1)^{i} p_i  y_i P_{jq}^{(i)}(y) ]
I_{\phi_q,\gamma(y)}(y).\leqno(2.11) $$ It is easily seen from the
following evident equality in view of the definition (2.10),
$$\int_{\partial \gamma(y)}\frac  {df_0 \wedge \cdots \wedge
df_{K-1}}{(f_0 -y_0)  \cdots (f_{K-1} - y_{K-1})} \wedge df_i \wedge
i_E(\psi_q)  =0,  \;0 \leq i \leq {K-1}.$$  $$\int_{\partial \gamma(y)}\frac
 {df_0 \wedge \cdots \wedge df_{K-1}}{(f_0 -y_0)  \cdots (f_{K-1} - y_{K-1})}
\wedge d\varphi  =0, \; \varphi \in  \Omega_X^{N-1}.$$  Let us compare the
relation
\begin{flushleft} $(2.12)\;\;\;\; d\int_{\gamma(y)}{\omega}_j =
\ell_j  \sum_{q=1}^{\mu(X_0)} [\sum_{i=0}^{K-1} (-1)^{i} P_{jq}^{(i)}(y)dy_i
] I_{\phi_q,\gamma(y)}(y),$ \end{flushleft}
obtained from (2.7) and (2.11). As
a result  we get equations between $I_{\phi_q}(y)$ and
$\frac{\partial}{\partial y_\ell} I_{\phi_q},$  $0 \leq \ell  \leq {K
-1},$ (we omit to specify $\gamma(y)$ except necessary
cases), \begin{flushleft}   $ \frac{\partial}{\partial y_\ell
}[\sum_{q=1}^{\mu(X_0)}  \sum_{i=0}^{K-1} (-1)^{i}p_i y_i P_{jq}^{(i)}
I_{\phi_q }]$  $= \ell_j  \sum_{q=1}^{\mu(X_0)} [\sum_{i=0}^{K-1} (-1)^{i}
P_{jq}^{(i)}(y)dy_i ] I_{\phi_q}(y), 1 \leq  j  \leq  \mu(X_0).
$ \end{flushleft} Thus we have obtained a system of differential equations to
be understood as the Gauss-Manin connexion of the singularity $X_0.$
To state the theorem in a simple form, we introduce
the following notations: $ {\bf I}_V = $
$(\int_{\gamma(y)}{\omega}_1,$  $ \cdots,
\int_{\gamma(y)}{\omega}_{\mu(X_0)}),$ $ {\bf I}_{\Phi} =$
 $(I_{\phi_1,\gamma(y)}(y), \cdots,$ $I_{\phi_{\mu(X_0)},\gamma(y)}(y) )
,$ i.e. vectors of integrals taken along a certain vanishing cycle
 $\gamma(y).$ We define several other $\mu(X_0) \times \mu(X_0)$ matrices  as
follows, $$ L_V = \diag(\ell_1, \cdots, \ell_{\mu(X_0)})$$ with  $\ell_i
=  w ({\omega}_i),$ $P^{(0)}(y) = (P^{(0)}_{jq}(y)), \cdots,$
$P^{({K-1})}(y) =
(P^{({K-1})}_{jq}(y)),1\leq j,q  \leq \mu(X_0),$  matrices consisting of
elements defined in (2.6).   \par In summing up the above arguments and the
theory due to Greuel \cite{Gr}, we obtain the following.
\begin{theorem}
 1).  For a quasihomogeneous mapping  $$ f: X \rightarrow Y$$  with isolated
complete intersection singularities  of dimension $N$ like (2.1), the
 Gauss-Manin system satisfied by  ${\bf I}_{\Phi}$  is described as
follows: \begin{flushleft}$(2.13)\;\;\;\; d[\sum_{i=0}^{K-1} (-1)^{i} p_i y_i
P^{(i)}(y){\bf I}_{\Phi}] =  L_V[\sum_{i=0}^{K-1} (-1)^{i} P^{(i)}(y)dy_i]
 {\bf I}_{\Phi}.$\end{flushleft}

\par 2).  The critical value $D$ (singular locus
of the system (2.13)) of   deformation  $X_y$ is given by   $D=\{y \in Y:
\Delta(y) =0 \}$ where   $$ \Delta(y) = \det (\sum_{i=0}^{K-1} (-1)^{i} p_i y_i
P^{(i)}(y)).\leqno(2.14)  $$
\par 3). The system of differential equations (2.13) is  a  holonomic
system.
\label{theorem21}  \end{theorem}     \par Let us
return to the problem (C.P.) of \S 1. Our main concern is to  understand the
integral (1.8) as a sum
of integrals like  (2.10) for  certain mapping $f.$ To adapt our (C.P.) to the
scheme explained before Theorem ~\ref{theorem21}, we treat the  following
mapping  $f: X \rightarrow Y$ for $ X = ( {\bf C}^{n+\mu}_u, 0), Y = ( {\bf
C}^{\mu+1}_y, 0).$ Concretely, it is defined as follows,  \begin{flushleft}
 $$\leqno(2.15) \left\{  \begin{array}{ccccc}  f_0(u) & = & F(z)&=  & y_0\\
 f_1(u) & = & z_{n+\mu}^P + { \langle z, \grad_z F(z) \rangle }^m  +
\sum_{i=1}^{\mu -1} z_{n+i}z^{\alpha^{(i)}} & = &y_1 \\  f_2(u) & = &
z_{n+1}&=  & y_2\\  \vdots & \vdots & \vdots&\vdots  & \vdots \\  f_{i+1}(u) &
= & z_{n+i}&=  & y_{i+1}\\  \vdots & \vdots & \vdots&\vdots  & \vdots \\
 f_\mu(u) & = & z_{n+\mu-1}&=  & y_\mu \\  \end{array}\right. $$
\end{flushleft} with notation $z=(z_1, \cdots, z_n), z'=(z_{n+1}, \cdots,
z_{n+\mu-1}), u = (z,z',z_{n+\mu}).$ Here the power $P$ is an integer that
corresponds to the denominator of $q \in \bf Q.$   \begin{lem}  For  $F(z)$
under conditions (C.2), (C.3), the mapping (2.15) $(f_0 ,$ $ \cdots ,$
$f_{\mu} )$ defines  an isolated quasihomogeneous complete intersection
singularity, $$X_0  =\{ u \in X : f_0(u) = \ldots = f_{\mu}(u) = 0  \}.$$
Namely
\begin{flushleft} $V  =$ $ \frac{\Omega^{n+\mu}_X} {
f_0\Omega^{n+\mu}_X +   f_1\Omega^{n+\mu}_X +  \sum_{i=1}^\mu
z_{i+n}\Omega^{n+\mu}_X  + dF \wedge \Omega^{n+\mu-1}_X + df_1\wedge
\Omega^{n+\mu-1}_X + \sum_{i=1}^\mu dz_{i+n} \wedge \Omega^{n+\mu-1}_X
},$ \end{flushleft}  is a finite dimensional vector space. \label
{lem22} \end{lem} {\bf Proof} The complete intersection property  follows from
the fact  that two polynomials $F(z)$ and $\langle z, \grad_z F(z) \rangle $
are of the same quasihomogeneous weight but with different coefficients. This
is a consequence of  (C.2), 1) which supposes that  $F(z)$ is not a
homogeneous polynomial.  The condition (C.3) entails immediately the
finite dimensionality of $V.$ {\bf Q.E.D.}
\par  To see that the components of
${\bf I}_{\Phi} $ defined for the mapping (2.15) give rise to
integrals of type (1.8), we prepare the following.
 \begin{lem}   1). Let us denote
$$
\Phi(z,z')  =
\left.\frac{\Omega^{n+\mu+1}_X} { df_0\wedge \ldots   df_\mu \wedge\Omega^{n}_X +
\sum_{i=0}^\mu f_{i}\Omega^{n+\mu+1}_X   }
\right|_{z_{n+\mu}=0}.
$$
Then
the following natural isomorphism holds $$
\Phi \cong \Phi(z,z') \otimes ( {\bf C}[z_{n+\mu}]/ <z_{n+\mu}^P>).$$
2). For  $\partial \gamma_{n-2}$ $\in H_{n+\mu -1}({\bf  C}^{n+\mu-1}
\setminus \cup_{i=0}^{\mu} \{ f_i = y_i \} \mid _{z_{n+\mu}=0
}, {\bf Z})$ a Leray coboundary of a vanishing cycle $\gamma_{n-2} \in
H_{n-2}(X_y\mid_{ z_{n+\mu}=0}, {\bf Z} )$  one can choose a
corresponding vanishing cycle ${\tilde \gamma}_{n-1} \in H_{n-1} ( X_y,
{\bf Z} )$ such that an equality $$ \leqno(2.16) \int_{\partial
\gamma_{n-2}} \phi(z) (f_1(z,y_2,\cdots, y_\mu,0)-y_1
)^{\frac{r+1}{P}-1} \frac{dz}{dF} = $$ $$=\epsilon (\frac{1}{2\pi
i})^\mu\int_{\partial\tilde  \gamma_{n-1}} \phi(z) z_{n+\mu}^r
\frac{du}{(f_0 -u_0) \cdots (f_{\mu}- u_\mu)}$$ holds, where $\epsilon
\in {\bf C}^{\times}$ such that $\epsilon^P =1.$ Furthermore, the cycle
$\partial\tilde  \gamma_{n-1}$ $\in H_{n+\mu}( {\bf
C}^{n+\mu}\setminus \cup_{i=0}^{\mu} \{f_i = u_i \}, {\bf Z})$ is
topologically equivalent to a product of a small circle on complex
$z_{n+\mu}-$ plane and $\partial \gamma_{n-2}.$ \label{lem23} \end{lem}
{\bf Proof} The statement 1) is evident. The statement 2) is an
integral version of statement 1), which can be shown by means of
 equality (2.10). {\bf Q.E.D.} \par Thus the singular locus of the
integral (1.8) can be given by that of $$\int_{\tilde\gamma_{n-1}}
\phi(z) z_{n+\mu}^r \frac{du}{df_0 \wedge \cdots \wedge df_{\mu}}$$
with ${\tilde \gamma}_{n-1}\in H_{n-1}(X_y,{\bf Z})$ after substitution
$y_1 = -W_1(x,t)$ (Case 1 after Lemma~\ref{lem12}), or $y_1 = 0$ (Case
2 after Lemma~\ref{lem12}),   $y_i = W_i(x,t), 2 \leq i \leq \mu.$ Let
us remind that we denoted the quasihomogeneous weight of function $f_i$
by $p_i,$ $0\leq i \leq \mu.$ We define matrices $P^{(i)}(y), 2 \leq i
\leq \mu$ for the mapping (2.15) after the master (2.6) and
Theorem~\ref{theorem21}. Combining Theorem~\ref{theorem21} with
Lemma~\ref{lem23}, we obtain the following.  \begin{theorem}  The
defining equation of wave front (Propositi on~\ref{prop11}, c)) is
given by the following polynomial: $$ \phi(x,t,s) = \det
(\sum_{i=0}^{\mu} (-1)^{i}  p_i y_i P^{(i)}(y))\mid_{y_0 =s, y_i =
W_i(x,t), 2 \leq i \leq \mu}. \leqno(2.17)  $$ Here the restriction shall be
imposed in accordance with   two cases treated just after Lemma~\ref{lem12}
i.e.  $y_1 = - W_1(x,t)$  in Case 1 and $ y_1 = 0$  in Case 2.
 \label{theorem22}  \end{theorem}

 {
\noindent
  \begin{flushleft}
\begin{minipage}[t]{6.2cm}   \begin{center}
{\footnotesize Moscow Independent University
,\\ Bol'shoj Vlasievskij pereulok 11,
121002, Moscow, Russia\\ E-mail:
tanabe@mpim-bonn.mpg.de,
tanabe@mccme.ru} \end{center} \end{minipage}\hfill
\end{flushleft}
\end{document}